\documentstyle[12pt]{article}

\begin{document}
\textheight 220mm \textwidth  140mm \oddsidemargin=0.05in
\evensidemargin=0.05in
\title{\bf  A partial order principle and vector variational principle for $\epsilon$-efficient solutions in the sense of N\'{e}meth$^{\footnotesize 1}$} \setcounter{footnote}{1} \footnotetext{ This
work was supported by the National Natural Science Foundation of
China (No. 11471236, 11561049).

{{\sl\ E-mail address:} qjhsd@sina.com, jhqiu@suda.edu.cn}}
\author{{Jing-Hui Qiu}\\
{\footnotesize  School of Mathematical Sciences, Soochow University,
Suzhou 215006,  PR China}
 }
\date{}
\maketitle

\baselineskip 16pt
\begin{quote}

{\bf Abstract.} \
In this paper, we establish a partial order principle, which is useful to deriving vector Ekeland variational principle (denoted by EVP).
By using the partial order principle and extending Gerstewitz's functions, we obtain a vector EVP for  $\epsilon$-efficient solutions in the sense of N\'{e}meth, which essentially  improves the earlier results by removing a usual assumption for boundedness of range of the objective function. From this, we also deduce several special vector EVPs, which improve and generalize the related known results.\\



\end{quote}

\section*{ \large\bf 1.  Introduction.   }

\hspace*{\parindent}  In 1972, Ekeland (see [13, 14]) presented
a variational principle, now known as Ekeland variational principle
(briefly, denoted by EVP), which says that for any lower
semi-continuous function $f$ bounded from below on a complete metric
space, there exists a slightly perturbed version of this function
that has a strict minimum. In the last four decades, the famous EVP
emerged as one of the most important results of nonlinear analysis
and its application covers numerous areas such as  optimization, optimal control theory,  fixed point theory,
nonsmooth analysis, Banach space geometry, game theory, nonlinear equations, dynamical systems, etc.; for example, see [3, 10, 14, 15, 19, 34, 49].  Motivated by
its wide usefulness, many authors have been interested in
extending EVP to the case with vector-valued maps or set-valued
maps; see, for example [2, 4-7, 9-12, 16, 17, 19, 20, 23-26, 28, 29, 32, 33, 38-40, 42-44, 46, 47, 50] and the references therein.

In this paper, we consider extensions of EVP when the objective function is a vector-valued map $f: (X, d) \rightarrow Y$, where $(X, d)$ is a complete metric space and $Y$ is a real
quasi-ordered (topological) vector space. A systematization of such results can be found in, for example, [10, 19, 20].
The common feature of these results is the presence of a certain term $d(x, x^{\prime})\,k_0$ in the perturbation, where $k_0\in D\backslash\{0\}$ and $D$ is an ordering cone.
Bednarczuk and Zagrodny (see [7, Theorem 4.1]) proved a vector EVP, where the perturbation is given by a bounded convex subset $H$ of the ordering cone $D$ multiplied by the distance function $d(x, x^{\prime})$, i.e., its  form is as  $d(x, x^{\prime})\, H$. This generalizes the case where directions of the perturbations are singleton $\{k_0\}$.
Tammer and Z$\breve{a}$linescu also considered this type of EVPs and gave an improvement of the above result; see [47, Theorem 6.2].
More generally, Guti\'{e}rrez, Jim\'{e}nez and Novo [23] introduced a set-valued metric, which takes values in the set family of all subsets  of the ordering cone and satisfies the triangle inequality.
By using it they gave an original approach to extending the scalar-valued EVP to a vector-valued map, where the perturbation contains a set-valued metric. From this, they deduced several special versions of EVP involving approximate solutions for vector optimization problems. However, in their work the assumption that the ordering cone $D$ is w-normal is required (see [23]). This requirement restricts the applicable extent of the new version of EVP.  Qiu [40] introduced a slightly more general notion: set-valued quasi-metrics, and introduced the notion of compatibility between a set-valued quasi-metric and the original metric $d$. By means of the notions,  Qiu proved a general  vector EVP, where the perturbation contains a set-valued quasi-metric compatible with  the original metric. Here, one needs not assume that the ordering cone is w-normal.
From the general EVP, Qiu deduced a number of special  vector EVPs, which improve the related known results.  Particularly, Qiu obtained several EVPs for $\epsilon$-efficient solutions in the sense of N\'{e}meth, which improve the related results in [23].

In order to express our purpose clearly, we recall some details on this topic. Let $Y$ be a locally convex Hausdorff topological vector space (briefly, denoted by locally convex space) and $Y^*$ be its topological dual.  For any $\xi\in
Y^*$, we define a continuous semi-norm $p_{\xi}$ on $Y$ as follows:
$p_{\xi}(y):=|\xi(y)|, \ \,\forall y\in Y.$ The semi-norm family
$\{p_{\xi}:\,\xi\in Y^*\}$ generates a locally convex Hausdorff
topology on $Y$ (see, e.g., [27, 30, 31, 48]), which is called the weak
topology on $Y$ and denoted by $\sigma(Y, Y^*)$. For any nonempty
subset $F$ of $Y^*$, the semi-norm family $\{p_{\xi}:\;\xi\in F\}$
can also generate a locally convex topology (which need not be
Hausdorff) on $Y$, which is denoted by $\sigma(Y, F)$. In [30], the
topology $\sigma(Y, F)$ is called the $F$-projective topology. If
$A,\, B\subset Y$ and $\alpha\in R$, the sets $A+B$ and $\alpha A$
are defined as follows:
$$A+B:=\,\{z\in Y:\, \exists x\in A,\,\exists y\in B \ {\rm such\ that}\  z=x+y\},$$
$$ \alpha A:=\{z\in
Y:\, \exists x\in A \ {\rm such\ that}\ z=\alpha x\}.$$
A nonempty subset $D$ of $Y$ is called
a cone if $\alpha D\subset D$ for any $\alpha\geq 0$. And $D$ is
called a convex cone if $D+D\subset D$ and $\alpha D\subset D$ for
any $\alpha\geq 0$. Moreover, a convex cone is called a pointed convex cone if $D\cap(-D) = \{0\}$.
A pointed convex cone $D$ can specify a partial order in $Y$ as follows.
$$y_1,\,y_2\in Y,\ \  y_1\leq_D y_2\ \Leftrightarrow\ y_1-y_2\in -D.$$
The positive polar cone of $D$  is denoted by $D^+$, that is,
$D^+\,=\,\{\xi\in Y^*:\, \xi(d)\geq 0, \ \forall d\in D\}$. For $H\subset D\backslash\{0\}$,  the set $\{\xi\in Y^*:\,\inf\{\xi(h): h\in H\} > 0\}$ is denoted by $H^{+s}$.
A nonempty set $M\subset Y$  is said to be $D$-bounded by scalarization (briefly, denoted by $D$-bounded) if (see [23, Definition 3.3])
$$\inf\{\xi(y):\, y\in M\}\, >\,-\infty,\ \ \forall \xi\in D^+.$$
 Let us consider the following
vector optimization problem:
$$   {\rm Min}\{f(x):\, x\in S\}, \eqno{(1)}$$
where $f: X\rightarrow Y$ is a vector-valued map and $S$ is a
nonempty closed subset of $X$. A point $x_0\in S$ is called an
efficient solution of (1) if
$$  (f(S)-f(x_0))\cap (-D\backslash\{0\})=\emptyset,$$
where $f(S)$ denotes the set $\cup_{x\in S}\{f(x)\}$.

Guti\'{e}rrez, Jim\'{e}nez and Novo  introduced the $(C, \epsilon)$-efficiency concept, which extends and unifies several $\epsilon$-efficiency notions (see[21, 22]).\\

{ DEFINITION 1.1} (see [22]).\ {\sl A nonempty set $C\subset Y$ is
coradiant if $\cup_{\beta \geq 1}\beta C = C$.}\\

{ DEFINITION 1.2} (see [22]).\ {\sl Let $D$ be an ordering cone,
$C\subset D\backslash\{0\}$ be a coradiant set and $\epsilon >0$. A
point $x_0\in S$ is a $(C,\epsilon)$-efficient solution of problem
{\rm (1)} if $(f(S)-f(x_0))\cap (-\epsilon C) =\emptyset.$
In this case, we also denote  $x_0\in AE(C, \epsilon)$
}\\

In particular, if $C:= H+D$, where $H\subset D\backslash\{0\}$, then we can easily verify that $C$ is a coradiant set and $C\subset D\backslash\{0\}$. Thus, we obtain the concept of approximate efficiency due to N\'{e}meth.\\

{ DEFINITION 1.3} (see [22, 23, 33, 40]).\ {\sl Let $H\subset
D\backslash\{0\}$ and $\epsilon >0$. A point $x_0\in S$ is said to
be an $\epsilon$-efficient solution of {\rm (1)} in the sense of
N\'{e}meth (with respect to $H$) if $(f(S)-f(x_0))\cap (-\epsilon H-D)\,=\,\emptyset$.
In this case, we also denote $x_0\in AE(C_H, \epsilon)$, where $C_H = H+D$.
}\\

Usually, we assume that $H\subset  D\backslash\{0\}$  is a $D$-convex set, i.e.,  $H+D$ is a convex set.  Let $H$ be a $D$-convex set and $\gamma >0$. For any  $x\in S$, put
$$S(x):=\,\{z\in S:\, f(x)\in f(z) +\gamma d(x,z) H +D\}.\eqno{(2)}$$
It is easy to verify that $x\in S(x)$ and  $S(z)\subset S(x)$ for every $z\in S(x)$.\\

{ DEFINITION 1.4}  (see [40]).\ {\sl  Let $X$ be a metric space and let
$S(\cdot ):\, X\rightarrow 2^X\backslash\{\emptyset\}$ be a
set-valued map. The set-valued map $S(\cdot )$ is said to be
dynamically closed at $x\in X$ if $(x_n)\subset S(x),\
S(x_{n+1})\subset S(x_n)\subset S(x)$ for all $n$ and
$x_n\rightarrow \bar{x}$ then $\bar{x}\in S(x)$. In this case, we
also say that $S(x)$ is dynamically closed.}\\

We remark that a property similar to the above dynamical closedness, i.e., the so-called limiting monotonicity property, was  also introduced in [4, 5].
Let's recall  the following assumption (see [40]):

{\rm (Q3)} \ For any $x\in S(x_0)$, $S(x)$ is dynamically closed.

 Now we can relate a vector EVP in [40] for $\epsilon$-efficient solutions in the sense of N\'{e}meth, which generalizes [23, Theorem 5.11].\\

{ THEOREM 1.5} (see [40, Theorem  6.3]).\  {\sl Let $H\subset D\backslash\{0\}$ be a
$D$-convex set, let $0\not\in {\rm cl}(H+D)$,  and let assumption {\rm (Q3)}
be satisfied by considering $S(\cdot)$ determined by {\rm (2)}.  Let $x_0\in S$
be an $\epsilon$-efficient solution of {\rm (1)} in the sense of
N\'{e}meth with respect to $H$,  and assume that  the set
$(f(S)-f(x_0))\cap (-\epsilon ({\rm cone}(C_H)\backslash C_H))$ is
$D$-bounded, where $C_H = H+D$ and ${\rm cone}(C_H)$ denotes the cone generated by $C_H$.  Then, for any $\gamma >0$,  there exists $\hat{x} \in S$ such that

{\rm (a)} \  $f(x_0)\,\in\, f(\hat{x}) +\gamma d(x_0, \hat{x}) H +D$;

{\rm (b)} \  $d(x_0, \hat{x}) H\cap (\epsilon/\gamma)\,({\rm
cone}(C_H)\backslash C_H)\,\not=\,\emptyset;$

{\rm (c)} \ $\forall x\in S\backslash\{\hat{x}\}, \  f(\hat{x})\,\not\in\, f(x) +\gamma d(\hat{x}, x) H +D.$}\\

As we have seen, whether in [40, Theorem 6.3] or in [23, Theorem 5.11],  the assumption that  $(f(S)-f(x_0))\cap (-\epsilon ({\rm cone}(C_H)\backslash C_H))$ is
$D$-bounded is necessary. In fact, in the proofs of the above two theorems, we need to verify that assumption (A6), i.e., $x_0\in AE(C_H, \epsilon)$  and  $(f(S)-f(x_0))\cap (-\epsilon({\rm cone}(C_H)\backslash C_H))$ being $D$-bounded, is satisfied (for more details, see  [23, 40]).   Thus, the above assumption on $D$-boundedness
is indispensable. In this paper, we shall follow another way of deriving this sort of results.  First, we establish a partial order principle,
  which consists of a partial order set $(X, \preceq)$
 and an extended real-valued function $\eta$ which is monotone with respect to $\preceq$.
 The partial order principle states that there exists a strong minimal point dominated by any given point
 provided that the monotone function $\eta$  satisfies three general conditions. This is indeed a variant of [42, Theorem 2.1]. By using the partial order principle and extending the Gerstewitz's function,  we obtain a vector EVP for $\epsilon$-efficient solutions in the sense of N\'{e}meth, which essentially improves Theorem 1.5. To one's surprise, we find out that even though  the assumption  $(f(S)-f(x_0))\cap (-\epsilon ({\rm cone}(C_H)\backslash C_H))$ being $D$-bounded is completely  removed, the result of Theorem 1.5 remains true. From this, we also deduce several results, which improve [40, Theorem 6.5] and [23, Theorem 5.12]. Moreover, from the partial order principle, we obtain a vector EVP, where the perturbation contains a $\sigma$-convex set, which improves  [7, Theorem 4.1], [47, Theorem 6.2] and [40, Theorem 6.8] by relaxing the lower boundedness on ranges of objective functions.

 This paper is structured as follows.
In section 2, we establish a partial order principle, which is useful to deriving vector EVPs.  In section 3, we consider extending Gerstewitz's functions from a singleton $\{k_0\}$ to a set $H$ and discuss their properties. In section 4, by using the partial order principle and generalized Gerstewitz's functions we obtain a vector EVP for $\epsilon$-efficient solutions in the sense of N\'{e}meth, which improves the earlier results by removing a certain assumption on $D$-boundedness of the range of the objective function.  From this, we also deduce several interesting EVPs, which improve related known results. Finally, in section 5, still using the partial order principle, we obtain a vector EVP, where the perturbation contains a  $\sigma$-convex set (i.e., cs-complete bounded set; see [47]). The EVP improves several known EVPs by relaxing the lower boundedness for the range of the objective function.\\

\section*{\large\bf  2.  A partial order principle.}

 \hspace*{\parindent} In this section, we present a partial order principle, which is a useful tool of deriving EVPs. In fact, it is a variety of the pre-order principle in [42].

  Let $X$ be a nonempty set. As in [17], a binary
relation $\preceq$ on $X$ is called a pre-order if it satisfies the
transitive property; a quasi order if it satisfies the reflexive and
transitive properties; a partial order if it satisfies the
antisymmetric, reflexive and transitive properties. Let $(X,
\preceq)$ be a partial order set. An extended real-valued function
$\eta:\, (X, \preceq) \rightarrow {\bf R}\cup\{\pm \infty\}$ is called
monotone with respect to $\preceq$ if for any $x_1,\, x_2\in X$,
$$x_1\preceq x_2\ \ \Longrightarrow\ \ \eta(x_1)\leq \eta(x_2).$$
For any given $x_0\in X$, denote $S(x_0)$ the set $\{x\in X:\,
x\preceq x_0\}$. First we give a partial order  principle as follows.\\

{ THEOREM 2.1.} (refer to [42, Theorem  2.1]).\
{\sl Let $(X, \preceq)$ be a partial order  set,
$x_0\in X$ be given  and $\eta:\, (X,
\preceq) \rightarrow R\cup\{\pm\infty\}$ be an extended real-valued
function which is monotone with respect to $\preceq$.

Suppose that the following conditions are satisfied:

{\rm (A)}  \  $-\infty < \inf\{\eta(x):\, x\in S(x_0)\} <+\infty$;

{\rm (B)}  \  for any $x\in S(x_0)\backslash\{x_0\}$ with $-\infty<\eta(x)<+\infty$
and any $x^{\prime}\in S(x)\backslash\{x\}$, one has
$\eta(x)>\eta(x^{\prime})$;

{\rm (C)} \ for any sequence $(x_n)\subset S(x_0)$ with $x_n\in
S(x_{n-1})\backslash\{x_{n-1}\},\ \forall n$, such that  $\eta(x_n) -\inf\{\eta(x):\,
x\in S(x_{n-1})\} \rightarrow 0$ \  $(n\rightarrow \infty)$, there
exists $u\in X$ such that $u\in S(x_n),\ \forall n$.

Then there exists $\hat{x}\in X$ such that

{\rm (a)} \ $\hat{x}\in S(x_0)$;

{\rm (b)} \ $S(\hat{x}) =  \{\hat{x}\}$.}\\

{\sl Proof.}\ For brevity,  denote $\inf\{\eta(x):\, x\in
S(x_0)\}$ by $\inf \eta\circ S(x_0)$. By (A), we know that
$$-\infty <\inf \eta\circ S(x_0)<+\infty. \eqno{(3)}$$
If $S(x_0) = \{x_0\}$, then we may take $\hat{x}:= x_0$. Clearly, it satisfies (a) and (b).
If $S(x_0) \not= \{x_0\}$, then by (3) we may take $x_1\in S(x_0)\backslash\{x_0\}$ such that
$$\eta(x_1) \,<\, \inf \eta\circ S(x_0) +\frac{1}{2}. \eqno{(4)}$$
By the transitive property of $\preceq$, we have
$$ S(x_1)\subset S(x_0). \eqno{(5)}$$
If  $S(x_1) = \{x_1\}$, then we may take  $\hat{x}:= x_1$.  Clearly, it satisfies  (a) and (b).
If $S(x_1) \not= \{x_1\}$, then by (3), (4) and (5) we conclude that
$$-\infty <\inf \eta\circ S(x_1)<+\infty.$$
We may take $x_2\in S(x_1)\backslash\{x_1\}$ such that
$$\eta(x_2)\, <\, \inf\eta\circ S(x_1) +\frac{1}{2^2}.$$
In general, let $x_{n-1}\in X\  (n\geq 1)$  be given.
If $S(x_{n-1})  = \{x_{n-1}\}$, then we may take $\hat{x}:= x_{n-1}$. Clearly, it satisfies (a) and (b).
If  $S(x_{n-1})\not=\{x_{n-1}\}$, then we  conclude that
$$-\infty \,<\,\inf\,\eta\circ S(x_{n-1})\,<\, +\infty.$$
 We may take $x_n\in S(x_{n-1})\backslash\{x_{n-1}\}$ such that
$$\eta(x_n) < \inf\eta\circ S(x_{n-1}) +\frac{1}{2^n}.  \eqno{(6)}$$
Without loss of generality, we  assume that  $S(x_n)\not=\{x_n\}$ for every $n$.
Thus, we obtain a sequence
 $(x_n)\subset S(x_0)$ with $x_n\in
S(x_{n-1})\backslash \{x_{n-1}\},\ \forall n$, such that
$$\eta(x_n) - \inf \eta\circ S(x_{n-1}) \,<\, \frac{1}{2^n}\, \rightarrow\, 0\  \  (n\rightarrow\infty).$$
 By (C), there exists
$\hat{x}\in X$ such that $$\hat{x}\in S(x_n), \ \ \forall n.
\eqno{(7)}$$
Clearly, $\hat{x}\in S(x_0)$, that is, $\hat{x}$
satisfies (a). Next, we show that $\hat{x}$ satisfies (b),  that is,
$S(\hat{x}) = \{\hat{x}\}$.  If not, there exists
$\bar{x}\in S(\hat{x})$ and $\bar{x}\not=\hat{x}$. By (B),
$$\eta(\hat{x}) >\eta(\bar{x}). \eqno{(8)}$$
On the other hand, by $\bar{x}\in S(\hat{x})$ and (7) we have
$$\bar{x}\in S(x_n),  \ \ \forall n.\eqno{(9)}$$
Since $\eta$ is monotone with respect to $\preceq$, by (7), (6)
and (9) we have
\begin{eqnarray*}
\eta(\hat{x})\, \leq\, \eta(x_n) \,&<&\,\inf\eta\circ S(x_{n-1})
+\frac{1}{2^n}\\
&\leq&\,\eta(\bar{x}) +\frac{1}{2^n},\ \forall n.
\end{eqnarray*}
Letting $n\rightarrow\infty$, we have
$\eta(\hat{x}) \,\leq\, \eta(\bar{x})$, which contradicts (8).  \hfill\framebox[2mm]{}\\

\section*{\large\bf  3. Generalized Gerstewitz's functions and their properties. }

\hspace*{\parindent}  A useful approach for solving a vector problem is to reduce it to
a scalar problem.  Gerstewitz's functions introduced in [18] are often used as the
basis of the scalarization. In the framework of topological vector
spaces, Gerstewitz's functions generated by  closed convex
(solid) cones and their properties have been investigated thoroughly,
for example, see [10, 18, 19, 43, 44] and the references
therein. In  this section, we  consider Gerstewitz's functions and their generalizations in a more general framework.

In the following, we always assume that $Y$ is a real vector space.
For a nonempty subset $A\subset Y$,  the vector closure of $A$ is defined
as follows (refer to [1, 41, 43]):
$${\rm vcl}(A)\,=\,\{y\in Y:\,\exists v\in Y, \exists\lambda_n\geq
0, \lambda_n \rightarrow 0\ {\rm such\ that}\ y+\lambda_n v\,\in\,A,
\forall n\in N\}.$$  For any given $v_0\in Y$, we define
the $v_0$-vector closure (briefly, $v_0$-closure) of $A$ as follows:
$${\rm vcl}_{v_0}(A)\,=\,\{y\in Y:\,\exists\lambda_n\geq0,
\lambda_n\rightarrow 0\ {\rm such\ that}\ y+\lambda_n v_0\in A,
\forall n\in N\}.$$ Obviously,
$$A\subset {\rm vcl}_{v_0}(A) \subset \cup_{v\in Y}\, {\rm vcl}_v(A)
= {\rm vcl}(A).$$ All the above inclusions are proper. Moreover, if
$Y$ is a Hausdorff topological vector space (briefly, denoted by t.v.s.) and ${\rm cl}(A)$ denotes the
closure of $A$, then ${\rm vcl}(A)\subset {\rm cl}(A)$ and the
inclusion is also proper.
A subset $A$ of $Y$ is said to be $v_0$-closed if $A={\rm vcl}_{v_0}(A)$; to be vectorially closed if $A ={\rm vcl}(A)$;  to be (topologically) closed if $A = {\rm cl}(A)$.
In general, a nonempty set $A\subset Y$ need not be $v_0$-closed; a $v_0$-closed set need not be vectorially closed and a vectorially closed set need not be topologically closed (for details, see [43]).
Let $D\subset Y$ be a convex cone  and $\leq_D$ be a quasi-order determined by $D$. In this case,  $D$ is called the ordering cone or positive cone.
 We always assume that $D$ is nontrivial, i.e., $D\not=\{0\}$ and $D\not=Y$. Let  $k_0\in D\backslash -D$ be given.
 For any $y\in Y$, if there exists $t\in {\bf R}$ such that $y\in t k_0-D$, then for any    $t^{\prime}>t$, $y\in t^{\prime} k_0-D$. Thus, we can define a function $\xi_{k_0}:\, Y \rightarrow {\bf R}\cup\{ \pm \infty\}$  as follows: if there exists $t\in {\bf R}$ such that $y\in t k_0-D$, then define $\xi_{k_0}(y)\,=\,\inf\{t\in {\bf R}:\, y\in t k_0-D\}$; or else, define $\xi_{k_0}(y) = +\infty$. Such a function is called a Gerstewitz's function generated by $D$ and $k_0$.

 The following results concerning Gerstewitz's functions originate from [18,  19].\\

 { PROPOSITION 3.1} (see [43, Lemma 2.6]). \ {\sl There exists $z\in Y$ such that $\xi_{k_0}(z) =-\infty$ iff $k_0\in -{\rm vcl}(D)$.}\\

 { PROPOSITION 3.2} (see [19, 43, 44]). \ {\sl Let $D\subset Y$ be a convex cone and $k_0\in D\backslash-{\rm vcl}(D)$. Then, the Gerstewitz's function $\xi_{k_0}$ has the following properties:

{\rm (i)} \  $y_1\leq_D y_2\ \ \Longrightarrow\ \ \xi_{k_0}(y_1)\leq \xi_{k_0}(y_2),\ \ \forall y_1, y_2\in Y;$

{\rm (ii)} \  $\xi_{k_0}(\alpha y)\,=\,\alpha\,\xi_{k_0}(y),\ \  \forall y\in Y, \ \forall \alpha\geq 0;$

{\rm (iii)} \  $\xi_{k_0}(y_1+y_2)\,\leq\, \xi_{k_0}(y_1) + \xi_{k_0}(y_2),\ \ \forall y_1, y_2\in Y;$

Let $y\in Y$ and $r\in {\bf R}$. Then, we have:

{\rm (iv)} \  $\xi_{k_0}(y)  < r\  \Longleftrightarrow \  y\in  r k_0-{\rm vint}_{k_0}(D)$,
where ${\rm vint}_{k_0}(D) = (0,+\infty) k_0 +D;$

{\rm (v)} \  $\xi_{k_0}(y)  \leq r\  \Longleftrightarrow \  y\in rk_0 -{\rm vcl}_{k_0}(D);$

{\rm (vi)} \  $\xi_{k_0}(y) = r\  \Longleftrightarrow\  y\in rk_0 -({\rm vcl}_{k_0}(D)\backslash {\rm vint}_{k_0}(D));$

Particularly, $\xi_{k_0}(0) = 0,\ \  \xi_{k_0}(k_0) = 1$;

{\rm (vii)} \ $\xi_{k_0}(y)\geq r \ \Longleftrightarrow\  y\not\in r k_0-{\rm vint}_{k_0}(D);$

{\rm (viii)} \  $\xi_{k_0}(y) >r\ \Longleftrightarrow\ y\not\in r k_0 -{\rm vcl}_{k_0}(D);$

{\rm (ix)} \ $\xi_{k_0}(y+ r k_0)\,=\, \xi_{k_0}(y) + r.$}\\

As we have seen, the Gersterwitz's function $\xi_{k_0}$ plays an important role in deriving EVPs where perturbations contain a singleton $\{k_0\}$.
Now, we consider EVPs where the $\{k_0\}$ in perturbations is replaced by a subset $H$ of the ordering cone $D$; for example, see Theorem 1.5.
Thus, we need to extend the notion of Gerstewitz's functions.

Let $H\subset  D\backslash-D$ be a $D$-convex set. For any $y\in Y$, if there exists $t\in {\bf R}$ such that $y\in t\,H-D$, then for any $t^{\prime} > t$, $y\in t^{\prime}\,H-D$. Thus, we can define a function $\xi_H:\, Y\rightarrow {\bf R}\cup\{\pm\infty\}$ as follows: if there exists $t\in {\bf R}$ such that $y\in t\,H-D$, then define $\xi_H(y) = \inf\{t\in{\bf R}:\, y\in t\,H-D\}$; or else, define $\xi_H(y) = +\infty$.
We call such a function a generalized Gerstewitz's function generated by $D$ and $H$.   Next, we give some properties of generalized Gerstewitz's functions.\\

 { PROPOSITION 3.3} \ {\sl There exists $z\in Y$ such that $\xi_H(z) =-\infty$ iff $0\not\in {\rm vcl}(H+D)$.}\\

{\sl Proof.} \  Assume that there exists $z\in Y$ such that $\xi_H(z) = -\infty$. Then, for any $n\in {\bf N}$, $z\in -n\, H-D$. Thus, $z/n\,\in\,-H-D$.
Letting  $n\rightarrow\infty$, we have
$$ 0\,\in\, {\rm vcl}(-H-D)\,=\,-{\rm vcl}(H+D).$$
Hence, $0\in {\rm vcl}(H+D)$.

Conversely, assume that $0\in {\rm vcl}(H+D)$. Then, there exists $v\in Y$ and a sequence $(\lambda_n)$ with $\lambda_n \geq 0$ and $\lambda_n \rightarrow 0$ such that $\lambda_n v\in H+D$. Since $0\not\in H+D$, we have $\lambda_n >0,\ \forall n.$  Thus,
$$-v\,\in\, -\frac{1}{\lambda_n} H -D, \ \  \forall n. $$
Put $z=-v$. Then $\xi_H(z) = -\infty$.\hfill\framebox[2mm]{}\\

 { PROPOSITION 3.4} \ {\sl Let $D\subset Y$ be a convex cone and $H\subset D\backslash -D$ be a $D$-convex set such that $0\not\in {\rm vcl}(H+D)$. Then, the the generalized Gerstewitz's function $\xi_H$ has the following properties:

 {\rm (i)} \  $y_1\leq_D y_2\ \ \Longrightarrow\ \ \xi_{H}(y_1)\leq \xi_{H}(y_2),\ \ \forall y_1, y_2\in Y;$

{\rm (ii)}  \  $\xi_H(0)\,=\, 0$;

{\rm (iii)} \  $\xi_{H}(\alpha y)\,=\,\alpha\,\xi_{H}(y),\ \  \forall y\in Y, \ \forall \alpha\geq 0;$

{\rm (iv)} \  $\xi_{H}(y_1+y_2)\,\leq\, \xi_{H}(y_1) + \xi_{H}(y_2)$  if $\xi_H(y_1) < 0$ and $\xi_H(y_2) < 0$.
}\\

{\sl Proof.} \  (i) \  Without loss of generality, we may assume that $\xi_H(y_2) <+\infty$.  For any $\epsilon >0$,
$$y_2\,\in\,(\xi_H(y_2) +\epsilon) H-D.$$
Since $y_1\leq_D y_2$, we have $y_1\in y_2-D$.
Thus, $$y_1\in (\xi_H(y_2) +\epsilon)H-D-D = (\xi_H(y_2) + \epsilon) H -D.$$
Hence $\xi_H(y_1) \leq \xi_H(y_2) + \epsilon$, which leads to $\xi_H(y_1) \leq \xi_H(y_2)$.

(ii) \  Obviously, $ 0\,\in\, 0\cdot H-D$, so $\xi_H(0) \leq 0$. Assume that $\xi_H(0) < 0$. Then, there exists $\epsilon >0$ such that $\xi_H(0) +\epsilon\,<\, 0$.
Thus,
$$0\,\in\,(\xi_H(0)+\epsilon)H-D\,=\,-(\xi_H(0)+\epsilon)(-H-D).$$
Since $-(\xi_H(y) +\epsilon) >0$, we have
$0\in -H-D$ and $0\in H+D$, contradicting $0\not\in {\rm vcl}(H+D)$.

(iii) \  If  $\alpha =0$, then from (ii), $\xi_H(\alpha y) = \xi_H(0)  =0$. Also, $\alpha\cdot \xi_H(y) = 0\cdot \xi_H(y) = 0$.
Hence, $\xi_H(\alpha y) = \alpha\, \xi_H(y)$ holds for $\alpha =0$.

If $\alpha >0$ and $\xi_H(y) < +\infty$, then for any $\epsilon >0$,
$y\,\in\, (\xi_H(y) +\epsilon) H -D.$
Thus,
$\alpha\,y\,\in\, \alpha(\xi_H(y) +\epsilon) H -D.$
Hence, $\xi_H(\alpha y) \,\leq\,\alpha\, \xi_H(y) +\alpha\,\epsilon.$
Since $\epsilon >0$ may be arbitrary small, we have
$$\xi_H(\alpha y) \,\leq\, \alpha\, \xi_H(y).\eqno{(10)}$$
Also,
$$\xi_H(y)\,=\,\xi_H(\frac{1}{\alpha}\,\alpha y)\,\leq\,\frac{1}{\alpha}\,\xi_H(\alpha y).$$
From this,
$$\xi_H(\alpha y) \,\geq\,\alpha\,\xi_H(y).\eqno{(11)}$$
Combining (10) and (11), we have
$\xi_H(\alpha y)\,=\, \alpha\, \xi_H(y).$

If $\alpha >0$ and  $\xi_H(y) = +\infty$, then for any $t\in {\bf R}$, $y\not\in t\,H-D$.  Thus, for any $t\in {\bf R}$, $\alpha\, y\,\not\in t\,H-D$.
 Hence $\xi_H(\alpha y) = +\infty$ and $\xi_H(\alpha y)\,=\,\alpha\, \xi_H(y)$ holds.

(iv) \  Assume that $\xi_H(y_1)  < 0$ and $\xi_H(y_2) <0$. Then, there exists $\epsilon >0$ such that $\xi_H(y_1) +\epsilon <0$ and $\xi_H(y_2) + \epsilon < 0$. Thus,
$$y_1\,\in\,(\xi_H(y_1)+\epsilon) H- D\,=\, -(\xi_H(y_1)+\epsilon)\,(-H-D)$$
and
$$y_2\,\in\,(\xi_H(y_2)+\epsilon) H- D\,=\, -(\xi_H(y_2)+\epsilon)\,(-H-D).$$
Since $-H-D$ is convex, we have
\begin{eqnarray*}
y_1 + y_2 \,&\in&\, -(\xi_H(y_1)+\epsilon)(-H-D)-(\xi_H(y_2) + \epsilon)(-H-D)\\
&=&\, -(\xi_H(y_1) +\xi_H(y_2) + 2\epsilon)(-H-D)\\
&=&\,(\xi_H(y_1) + \xi_H(y_2) + 2\epsilon) H-D.
\end{eqnarray*}
From this, $$\xi_H(y_1+y_2)\,\leq\, \xi_H(y_1) + \xi_H(y_2) + 2 \epsilon.$$
Since $2\epsilon >0$ may be arbitrary small, we have
$$ \xi_H(y_1+y_2) \,\leq\, \xi_H(y_1) + \xi_H(y_2).$$\hfill\framebox[2mm]{}\\

\section*{\large\bf  4. Vector EVPs  for $\epsilon$-efficient solutions in the sense of N\'{e}meth.}

\hspace*{\parindent} In this section, we assume that $(X,d)$ is a complete metric space, $Y$ is a  vector space quasi-ordered by a convex cone $D$, $f:\, X\rightarrow Y$ is a vector-valued map and $S$ is a nonempty closed subset of $X$. By using the partial order principle (i.e., Theorem 2.1) and generalized Gerstewitz's functions, we obtain the following vector EVP, which improves Theorem 1.5 by removing the assumption that the set $(f(S)- f(x_0)) \cap (-\epsilon({\rm cone}(C_H)\backslash C_H))$ is $D$-bounded.  Besides, the assumption in Theorem 1.5 that $0\not\in {\rm cl}(H+D)$ is replaed by a weaker one that $0\not\in {\rm vcl}(H+D)$.\\

{ THEOREM 4.1.}   \  {\sl Let $H\subset D$ be a
$D$-convex set such that $0\not\in {\rm vcl}(H+D)$, and let  assumption {\rm (Q3)}
be satisfied by considering $S(\cdot)$ determined by {\rm (2)}.  Let $x_0\in S$
be an $\epsilon$-efficient solution of {\rm (1)} in the sense of
N\'{e}meth with respect to $H$.  Then, for any $\gamma >0$,  there exists $\hat{x} \in S$ such that

{\rm (a)} \  $f(x_0)\,\in\, f(\hat{x}) +\gamma d(x_0, \hat{x}) H +D$;

{\rm (b)} \  $d(x_0, \hat{x}) H\cap (\epsilon/\gamma)\,({\rm
cone}(C_H)\backslash C_H)\,\not=\,\emptyset;$

{\rm (c)} \ $\forall x\in S\backslash\{\hat{x}\}, \  f(\hat{x})\,\not\in\, f(x) +\gamma d(\hat{x}, x) H +D.$}\\

{\sl Proof.} \  For $x, x^{\prime}\in S$, define $x^{\prime}\preceq x$ iff $f(x)\in f(x^{\prime}) + \gamma d(x, x^{\prime}) H +D.$
It is easy to verify that $\preceq$ is a partial order on $S$.  If $x^{\prime}\preceq x$ and $ x^{\prime}\not= x$, we denote $x^{\prime}\prec x$.
Define an extended real-valued function $\eta:\, (S, \preceq) \rightarrow {\bf R}\cup\{\pm\infty\}$ as follows:
$$\eta(x):=\,\xi_H(f(x)-y_0),\ \ x\in S,$$
where $y_0 = f(x_0)$.
Since $0\not\in {\rm vcl}(H+D)$, by Proposition 3.3, $\xi_H(f(x) - y_0)\not=-\infty$, that is, $\eta(x) \not= -\infty,\ \forall x\in S.$
Let $x^{\prime}\preceq x$. Then $$f(x) \in f(x^{\prime}) +\gamma d(x, x^{\prime}) H +D.$$
Thus,
$$f(x^{\prime})-f(x)\,\in\,-\gamma d(x, x^{\prime})H-D\,\subset\,-D$$
and
$$f(x^{\prime})-y_0\,\leq_D\, f(x)- y_0.$$
By Proposition 3.4(i), we have
$$\xi_H(f(x^{\prime})-y_0)\,\leq\,\xi_H(f(x)-y_0),\ \ {\rm that\ is, } \ \ \eta(x^{\prime}) \leq \eta(x).$$
Hence, $\eta$ is monotone with respect to $\preceq$.
We denote the set $\{x^{\prime}\in X:\, x^{\prime}\preceq x\}$ by $S(x)$.
Next, we prove that assumptions (A), (B) and (C) in Theorem 2.1 are satisfied.

{\bf Prove that (A) is satisfied.} \  Since $y_0= f(x_0)\not\in f(S) +\epsilon H +D$,
for any $x\in S(x_0)\subset S$,
$$f(x) - y_0 \,\not\in\,-\epsilon H-D,\ \ {\rm so}  \ \  \eta(x) = \xi_H(f(x)- y_0) \geq -\epsilon.$$
Also, by Proposition 3.4(ii),
$$\eta(x_0)\,=\, \xi_H(f(x_0)-y_0)\,=\,\xi_H(0)\,=\,0.$$
Hence,
$$-\infty < -\epsilon \leq \inf\{\eta(x):\, x\in S(x_0)\} \leq  \eta(x_0) = 0 <+\infty.$$
That is, (A) is satisfied.

{\bf Prove that (B) is satisfied.} \ Let $x\in S(x_0)\backslash\{x_0\}$ with $-\infty < \eta(x) <+\infty$ and let $x^{\prime}\in S(x) \backslash\{x\}$.
By  $x\in S(x_0)\backslash\{x_0\}$, we have
$$f(x_0) \,\in\, f(x) +\gamma d(x_0, x) H +D\ \  {\rm and} \ \  x\not= x_0.\eqno{(12)}$$
By  $x^{\prime}\in S(x)\backslash\{x\}$, we have
$$f(x)\,\in\, f(x^{\prime}) +\gamma d(x, x^{\prime}) H +D \ \ {\rm and} \ \ x^{\prime}\not= x. \eqno{(13)}$$
By (12), we know that
$$f(x) - y_0\,\in\, -\gamma d(x_0, x) H-D$$
and so
$$\xi_H(f(x)-y_0)\,\leq\,-\gamma d(x_0, x)  < 0. \eqno{(14)}$$
By (13), we know that
$$f(x^{\prime}) -  f(x)\,\in\, -\gamma d(x, x^{\prime}) H-D$$
and so
$$\xi_H(f(x^{\prime})- f(x))\,\leq\,-\gamma d(x, x^{\prime})  < 0. \eqno{(15)}$$
Remarking (14) and (15), and using Proposition 3.4(iv), we have
$$\xi_H(f(x^{\prime})-y_0) \,\leq\, \xi_H(f(x^{\prime})- f(x)) + \xi_H(f(x)- y_0).$$
From this and using (15), we have
$$\xi_H(f(x^{\prime})-y_0) - \xi_H(f(x) -y_0)\,\leq\, \xi_H(f(x^{\prime})-f(x))\,\leq\, -\gamma d(x, x^{\prime}).$$
That is,
$$\eta(x^{\prime}) - \eta(x)\,\leq\, -\gamma d(x, x^{\prime})$$ and so
$$\eta(x^{\prime}) \,\leq\, \eta(x) -\gamma d(x, x^{\prime})\,<\,\eta(x).$$
Thus, (B) is satisfied.

{\bf Prove that (C) is satisfied.} \  Let a sequence $(x_n) \subset S(x_0)$ with $x_n\in S(x_{n-1})\backslash\{x_{n-1}\}, \ \forall n,$ such that $\eta(x_n) - \inf\{\eta(x):\, x\in S(x_{n-1})\} \rightarrow 0\ (n\rightarrow\infty)$.
Since $x_0\succ x_1\succ x_2\succ \cdots \succ x_n\succ\cdots$, by the transitive property and the antisymmetric property, we have $x_m\prec  x_n,\ \forall m>n.$
That is,
$$f(x_n)\,\in\, f(x_m) + \gamma d(x_n, x_m) H +D.$$
From this,
$$\xi_H(f(x_m) -f(x_n))\,\leq\, -\gamma d(x_n, x_m)$$
and so
$$\gamma d(x_n, x_m) \,\leq\, -\xi_H(f(x_m)-f(x_n)).\eqno{(16)}$$
Since $\xi_H(f(x_n)- y_0) \leq -\gamma d(x_0, x_n) < 0$ and $ \xi_H(f(x_m)- f(x_n)) \leq -\gamma d(x_n, x_m) <0$,
by Proposition  3.4(iv), we have
$$\xi_H(f(x_m)-y_0)\,\leq\,\xi_H(f(x_m)-f(x_n)) + \xi_H(f(x_n) -y_0)$$
and so
$$-\xi_H(f(x_m)-f(x_n))\,\leq\,\xi_H(f(x_n)-y_0)-\xi_H(f(x_m)-y_0).\eqno{(17)}$$
Combining (16) and (17), and remarking that $x_m\in S(x_{n-1})$, we have
\begin{eqnarray*}
\gamma d(x_n, x_m)\,&\leq&\,\xi_H(f(x_n) -y_0) - \xi_H(f(x_m) - y_0)\\
&=&\,\eta(x_n) -\eta(x_m)\\
&\leq&\, \eta(x_n) -\inf \eta\circ S(x_{n-1})\ \rightarrow\ 0\ \  \  (n\rightarrow\infty).
\end{eqnarray*}
Hence $(x_n)$ is a Cauchy sequence.  Since $(X, d)$ is complete and $S\subset X$ is closed, there exists $\hat{x}\in S$ such that $x_n \rightarrow \hat{x} \ (n\rightarrow\infty).$
For any given $n$, $S(x_n)\subset S(x_0)$. we observe  that $(x_{n+p})_{p\in {\bf N}}\subset S(x_n)$ and $x_{n+p+1}\in S(x_{n+p})  \  \forall p.$
Since $x_{n+p} \rightarrow \hat{x} \ (p\rightarrow\infty)$ and $S(x_n)$ is dynamically closed by (Q3), we have $\hat{x}\in S(x_n),\ \forall n.$ That is, (C) is satisfied.

Now, applying Theorem 2.1 we conclude that there exists $\hat{x}\in S$ such that $\hat{x} \in S(x_0)$ and $S(\hat{x}) = \{\hat{x}\}$.
That is, $\hat{x}$ satisfies (a) and (c).  Finally we show that $\hat{x}$ satisfies (b). By (a),
$$f(x_0)\,\in\, f(\hat{x}) + \gamma d(x_0, \hat{x}) H +D.$$
Hence, there exists $h_0\in H$ and $d_0\in D$ such that
$$f(x_0)\,=\, f(\hat{x}) +\gamma d(x_0, \hat{x}) h_0 +d_0. \eqno{(18)}$$
Clearly,
$$d(x_0, \hat{x}) h_0\,\in\, d(x_0, \hat{x}) H \eqno{(19)}$$
and
$$d(x_0, \hat{x}) h_0\,\in\, {\rm cone}(H+D).\eqno{(20)}$$
Next we show that
$$d(\emph{}x_0, \hat{x}) h_0\,\not\in\,\frac{\epsilon}{\gamma}(H+D).\eqno{(21)}$$
Assume that $$d(x_0, \hat{x}) h_0\,\in\, \frac{\epsilon}{\gamma} (H+D).$$
Then
$$\gamma d(x_0, \hat{x}) h_0\,\in\, \epsilon (H+D) \,=\, \epsilon H +D.$$
Thus,
$$\gamma d(x_0, \hat{x}) h_0+ d_0\,\in\,  \epsilon H +D+D =\epsilon H +D.$$
Combining this with (18), we have
$$f(x_0) - f(\hat{x})\,\in\, \epsilon H +D,$$
which contradicts the assumption that
$$f(x_0)\,\not\in\, f(S) + \epsilon H +D.$$
Now, combining (19), (20) and (21), we have
\begin{eqnarray*}
d(x_0, \hat{x}) h_0\,&\in&\, d(x_0, \hat{x})H\cap ({\rm cone}(H+D)\backslash\frac{\epsilon}{\gamma}(H+D))\\
&=&\, d(x_0, \hat{x}) H \cap (\epsilon/\gamma)({\rm cone}(H+D)\backslash(H+D)).
\end{eqnarray*}
This means that (b) is satisfied.\hfill\framebox[2mm]{}\\

As in [23, 40], a vector-valued map $f:\, X\rightarrow Y$ is said to be sequentially submonotone with respect to $D$ (briefly, denoted by submonotone) if for every $x\in X$ and for each sequence $(x_n)$ such that $x_n \rightarrow x$ and
$f(x_m) \leq_D f(x_n),\ \forall m>n,$ it follows that $f(x)\leq_D f(x_n), \ \forall n.$
Sometimes, a submonotone vector-valued map is said to be $D$-sequentially lower monotone (briefly, denoted by $D$-slm or slm);  see, for example [25]. In [7], a submonotone vector-valued map is called a monotonically semi-continuous (denoted by msc) with respect to $D$ map; in [20] it is called a map with property (H4); and in [29] it is called a lower semi-continuous from above (briefly, denoted by lsca).
Let us observe that a $D$-lower semi-continuous vector-valued map $f:\, X\rightarrow Y$, i.e., $f$ such that the sets $\{x\in X:\, f(x)\leq_D y\}$ are closed for all $y\in Y$, is submonotone. But the converse is not true
even $Y$  is the real number space with the usual order, for example, see [8].

Next, we present  a particular version of vector EVP for $\epsilon$-efficient  solutions by giving a certain condition for (Q3) fulfilled.\\

{ THEOREM 4.2.}   \  {\sl  Let $H\subset D$ be a $D$-convex set such that $0\not\in{\rm vcl}(H+D)$, and let $x_0\in S$ be an $\epsilon$-efficient solution of {\rm (1)} in the sense of N\'{e}meth
with respect to $H$. Moreover, assume that $H+D$ is $h_0$-closed for a certain $h_0\in H$ and $f$ is submonotone.
Then, for any $\gamma >0$,  there exists $\hat{x} \in S$ such that

{\rm (a)} \  $f(x_0)\,\in\, f(\hat{x}) +\gamma d(x_0, \hat{x}) H +D$;

{\rm (b)} \   $d(x_0, \hat{x})\,<\, \epsilon/\gamma;$

{\rm (c)} \ $\forall x\in S\backslash\{\hat{x}\}, \  f(\hat{x})\,\not\in\, f(x) +\gamma d(\hat{x}, x) H +D.$}\\

{\sl Proof.} \  From Theorem 4.1, we only need to prove that (Q3) is satisfied.
Let $x\in S(x_0)$, $(x_n) \subset S(x)$ with $x_{n+1}\in S(x_n)$ and $x_n \rightarrow u$.
For any given $n$ and for every $m> n$, we have $x_m\in S(x_n)$ and hence $f(x_m) \leq_D f(x_n)$.  Since $f$ is submonotone and $x_m\rightarrow u\ (m\rightarrow\infty)$, we have $f(u) \leq_D f(x_n)$.
For $m>n$, $x_m\in S(x_n)$. Thus,
\begin{eqnarray*}
f(x_n)\,&\in&\, f(x_m) + \gamma d(x_n, x_m) H +D\\
&\subset&\, f(u) + \gamma d(x_n, x_m) H +D.
\end{eqnarray*}
Next, we show the result according to the following two cases.

Case 1. \  There exists $m>n$ such that $d(x_n, x_m) \geq d(x_n, u)$.
Then
\begin{eqnarray*}
f(x_n)\,&\in&\, f(u) + \gamma d(x_n, x_m) H +D\\
&\subset&\, f(u) + \gamma d(x_n, u) H +D.
\end{eqnarray*}
That is, $u\in S(x_n) \subset S(x)$.

Case 2. \  For every $m>n$,  $d(x_n, x_m) < d(x_n, u)$.
Then, from
$$f(x_n) \,\in\, f(u) +\gamma d(x_n, x_m) H +D,$$
we have
\begin{eqnarray*}
&\ &\, f(x_n) +\gamma(d(x_n, u) - d(x_n, x_m)) h_0\\
&\in&\, f(u)+\gamma d(x_n, x_m) H +D +\gamma(d(x_n, u) - d(x_n, x_m))h_0\\
&\subset &\,  f(u)+\gamma d(x_n, x_m) (H +D) +\gamma(d(x_n, u) - d(x_n, x_m))(H+D)\\
&=&\, f(u) + \gamma d(x_n, u) (H+D)\\
&=&\, f(u) + \gamma d(x_n, u) H +D.
\end{eqnarray*}
From this,
$$f(x_n) - f(u) +\gamma(d(x_n, u)-d(x_n, x_m)) h_0\,\in\, \gamma d(x_n, u) H + D.\eqno{(22)}$$
Since $d(x_n, u) - d(x_n, x_m) \rightarrow 0\ (m\rightarrow \infty)$ and  $\gamma d(x_n, u) H +D$ is $h_0$-closed,
by (22) we have
$$f(x_n) - f(u)\,\in\,\gamma d(x_n, u) H +D$$
and
$$f(x_n)\,\in\, f(u) +\gamma d(x_n, u) H +D.$$
That is, $u\in S(x_n) \subset S(x)$.
Thus, we have shown that (Q3) is satisfied.
Applying Theorem 4.1, there exists $\hat{x}\in S$ such that (a) and (c) are satisfied.
Next we show that  (b) is satisfied. If not, assume that $d(x_0, \hat{x}) \geq \epsilon/\gamma$. Then from (a), we have
\begin{eqnarray*}
f(x_0)\,&\in&\, f(\hat{x}) +\gamma d(x_0, \hat{x}) H +D\\
\,&\subset&\, f(\hat{x}) +\gamma\,\frac{\epsilon}{\gamma} H +D\\
&=&\, f(\hat{x}) + \epsilon H +D,
\end{eqnarray*}
which contradicts the assumption that $x_0\in S$ is an $\epsilon$-efficient solution of {\rm (1)} in the sense of N\'{e}meth
with respect to $H$, i.e.,
$f(x_0)\,\not\in\, f(S) +\epsilon H +D.$\hfill\framebox[2mm]{}\\

{ THEOREM 4.3.}   \  {\sl  Let $Y$ be a locally convex space, $D\subset Y$ be a closed convex cone and   $H\subset D\backslash -D$ be a $\sigma(Y, D^+)$-countably compact, $D$-convex set.
 Suppose that $f:\, X\rightarrow Y$ is submonotone and $x_0$ is an $\epsilon$-efficient solution of {\rm (1)} in the sense of
  N\'{e}meth with respect to $H$.
Then, for any $\gamma >0$,  there exists $\hat{x} \in S$ such that

{\rm (a)} \  $f(x_0)\,\in\, f(\hat{x}) +\gamma d(x_0, \hat{x}) H +D$;

{\rm (b)} \   $d(x_0, \hat{x})\,<\, \epsilon/\gamma;$

{\rm (c)} \ $\forall x\in S\backslash\{\hat{x}\}, \  f(\hat{x})\,\not\in\, f(x) +\gamma d(\hat{x}, x) H +D.$}\\

{\sl Proof.} \  By Theorem 4.2, we only need to prove that $H+D$ is  vectorially closed.
Let $z\in {\rm vcl}(H+D)$. Then, there exists $v_0\in Y$ and a sequence $(\epsilon_n)$ with $\epsilon_n>0$ and $\epsilon_n \rightarrow 0$ such that
$z+\epsilon_n v_0\,\in\, H+D$.
For each $n$, there exists $h_n\in H$ such that
$$z+\epsilon_n v_0\,\in\, h_n+D,$$
that is,
$$z-h_n +\epsilon_n v_0\,\in\, D.\eqno{(23)}$$
Since $H$ is $\sigma(Y, D^+)$-countably compact, the sequence $(h_n)\subset H$  has a $\sigma(Y, D^+)$-cluster point $h^{\prime}\in H$.
Take  any $\xi\in D^+$. From (23), we have
$$\xi(z) -\xi(h_n) +\epsilon_n \xi(v_0) \,\geq\, 0.\eqno{(24)}$$
Since a continuous map preserves cluster points, $\xi(h^{\prime})$ is a cluster point of $(\xi(h_n))_n$ in ${\bf R}$. Hence, there exists a subsequence $n_1 < n_2 < n_3<\cdots$ such that $\xi(h_{n_i})\,\rightarrow\,\xi(h^{\prime})\ (i\rightarrow\infty).$
By (24), we have
$$\xi(z) - \xi(h_{n_i}) + \epsilon_{n_i}\xi(v_0)\,\geq\, 0.$$
Letting $i\rightarrow\infty$, we have
$$\xi(z) -\xi(h^{\prime})\,\geq\, 0,\ \ {\rm i.e.,} \ \  \xi(z-h^{\prime})\,\geq\, 0.$$
Since $\xi\in D^+$ is arbitrary and $D$ is a closed convex cone, we have
$$ z-h^{\prime} \,\in\, D^{++} = D.$$
That is,
$z\,\in\, h^{\prime} +D \subset H+D$.
Thus, we have shown that $H+D \,=\, {\rm vcl}(H+D)$. Now, from Theorem 4.2 we obtain the result.\hfill\framebox[2mm]{}\\

{ Remark 4.4.}  \   [40, Theorem 6.5]  also gives the same result as in Theorem 4.3, but there one needs to assume that $H$ is a base of $D$. Here, we have removed the assumption. Clearly, Theorem 4.3 improves [40, Theorem 6.5] and also improves  [23,  Theorem 5.12].\\

\section*{\large\bf  5. Vector EVP   with perturbation containing a $\sigma$-convex set.}

\hspace*{\parindent}  Vector EVPs, where perturbations are of type $d(x,y) H$, are also considered by Bednarczuk and Zagrodny [7], Tammer and Z$\breve{a}$linescu [47] and  Qiu [40].
For details, see [7, theorem 4.1], [47, Theorem 6.2] and [40, Theorem 6.8].  We shall see that our partial order principle, i.e.,  Theorem 2.1, also implies this type of EVPs.
We shall obtain a vector EVP, where the perturbation contains a $\sigma$-convex set, which improves the above three results.  First, we recall some  terms and notions.
Let $Y$ be a t.v.s. and $B\subset Y$ be nonempty.  A convex series of points of $B$ is a series of the form $\sum_{n=1}^{\infty} \lambda_n b_n$, where every $b_n\in B$, every $\lambda_n \geq 0$  and $\sum_{n=1}^{\infty}\lambda_n\,=\,1$. $B$  is said to be a $\sigma$-convex  if every convex series of its points converges to a point of $B$ (see [35, 42]).
In fact, we can easily prove that a set is $\sigma$-convex iff it is cs-complete and bounded (see [47, 49]). Let $B$ be a $\sigma$-convex set. Then, for a sequence $(b_n)$ in $B$ and a real sequence $(\lambda_n)$ with $\lambda_n\geq 0$ and $0 < \sum_{n=1}^{\infty} \lambda_n < +\infty$,
$\sum_{n=1}^{\infty} \lambda_n b_n/\sum_{n=1}^{\infty}\lambda_n$ is a convex series in $B$ and it converges to some point $\bar{b}\in B$.
Thus, $\sum_{n=1}^{\infty}\lambda_n b_n$ converges to $(\sum_{n=1}^{\infty}\lambda_n) \bar{b}\,\in\,(\sum_{n=1}^{\infty}\lambda_n) B$.
We call a set $B$ sequentially complete iff every Cauchy sequence  $(b_n)$ in $B$, converges to a point of $B$.
 In [7], `` sequentially complete" is called ``semi-complete". It is easy to show that every sequentially complete, bounded  convex set  is a $\sigma$-convex set (see  [47, Remark 6.1]).
 However, a $\sigma$-convex set need not be sequentially complete. For example, an open ball $B$ in a Banach space is $\sigma$-convex, but it is not closed and hence is not sequentially complete (for details,  see [35, 41]).\\

{THEOREM 5.1.}   \  {\sl Let $(X, d)$ be a complete metric space, $Y$ be a t.v.s., $D\subset Y$ be a closed convex cone, $H\subset D$ be a $\sigma$-convex set  such that $0\not\in{\rm vcl}(H+D)$ and let $f:\, X\rightarrow Y$ be a submonotone vector-valued map. Suppose that  $x_0\in X$ and  $\epsilon >0$  such that
$$(f(x_0) -\epsilon H-D) \cap  f(X)\,=\,\emptyset.$$
Then, for any $\gamma >0$,  there exists $\hat{x} \in S$ such that

{\rm (a)} \  $f(x_0)\,\in\, f(\hat{x}) +\gamma d(x_0, \hat{x}) H +D$;

{\rm (b)} \   $d(x_0, \hat{x})\,<\, \epsilon/\gamma;$

{\rm (c)} \ $\forall x\in S\backslash\{\hat{x}\}, \  f(\hat{x})\,\not\in\, f(x) +\gamma d(\hat{x}, x) H +D.$}\\

{\sl Proof.} \  For $x, x^{\prime} \in X$, define $x^{\prime}\preceq x$ iff $f(x) \,\in\, f(x^{\prime}) + \gamma d(x, x^{\prime}) H +D$.
It is easy to show that $\preceq$ is a partial order on $X$. Here, in order to show that $\preceq$ satisfies antisymmetric property we only need to assume that $0\not\in H+D$. Hence we need not assume that $D$ is pointed.
As in the proof of Theorem 4.1, we define an extended real-valued function $\eta:\, (X,\preceq) \rightarrow {\bf R}\cup\{\pm\infty\}$ as follows:
$$\eta(x) :=\, \xi_H(f(x) -y_0),
 \  \   x\in X,$$ where $y_0 = f(x_0)$. For any $x\in X$, put $S(x):=\,\{x^{\prime}\in X:\, x^{\prime}\preceq x\}$. It's easy to prove that $\eta$  is monotone  with respect to  $\preceq$, and assumptions  (A) and (B) are satisfied. It suffices to prove that assumption (C) is satisfied.  Let a sequence  $(x_n) \subset S(x_0)$ with  $x_n\in S(x_{n-1})\backslash\{x_{n-1}\}$, $\forall n$, such that
 $$\eta(x_n) - \inf\{\eta(x):\, x\in S(x_{n-1})\}\,\rightarrow\,0\ \ \   (n\rightarrow \infty).$$
By $x_i \in S(x_{i-1})$ for $i=1,2, \cdots, n,$ we have
\begin{eqnarray*}
f(x_0)\,&\in&\, f(x_1) + \gamma d(x_0, x_1) H +D,\\
f(x_1)\,&\in &\, f(x_2) + \gamma d(x_1, x_2) H +D,\\
\cdots\,&\  &\, \cdots \  \  \  \  \   \  \  \   \cdots  \  \  \  \   \   \  \   \   \cdots\\
f(x_{n-1})\,&\in&\, f(x_n) + \gamma d(x_{n-1}, x_n) H +D.
\end{eqnarray*}
By adding the two sides of the above $n$ belonging relations, we have
$$f(x_0)\,\in\, f(x_n)  + \gamma\left(\sum_{i=1}^n d(x_{i-1}, x_i)\right) H +D$$
and
$$f(x_n) - y_0 \,=\, f(x_n) - f(x_0)\,\in\, -  \gamma\left(\sum_{i=1}^n d(x_{i-1}, x_i)\right) H -D.$$
From this,
$$\xi_H(f(x_n) -y_0)\,\leq\,-\gamma \left(\sum_{i=1}^n d(x_{i-1}, x_i)\right).$$
Hence,
$$\gamma\left(\sum_{i=1}^n d(x_{i-1}, x_i)\right)\,\leq\, -\xi_H(f(x_n)-y_0).\eqno{(25)}$$
By the assumption,
$$y_0 = f(x_0) \,\not\in\, f(x_n) + \epsilon H +D,$$
so
$$f(x_n) -y_0\,\not\in\,-\epsilon H -D.$$
Thus,
$$\xi_H(f(x_n) - y_0)\,\geq\,-\epsilon$$
and
$$-\xi_H(f(x_n) - y_0)\,\leq\,\epsilon.\eqno{(26)}$$
Combining (25) and (26), we have
$$\gamma\left(\sum_{i=1}^n d(x_{i-1}, x_i)\right) \,\leq\, \epsilon$$
and $$\sum_{i=1}^n d(x_{i-1}, x_i)\,\leq\, \frac{\epsilon}{\gamma},\ \ \ \forall n.$$
Thus,
$$\sum_{i=1}^{\infty} d(x_{i-1}, x_i)\,\leq\, \frac{\epsilon}{\gamma}.$$
Since $(X, d)$ is complete, there exists $u\in X$ such that $x_n \rightarrow u$ in $(X, d)$.
Next, we show that $u\in S(x_n),\ \forall n.$
 By $x_{i+1} \in S(x_i)$ for $i=n,\, n+1,\, \cdots,\, n+k-1$, we have
\begin{eqnarray*}
f(x_n)\,&\in&\, f(x_{n+1}) + \gamma d(x_n, x_{n+1}) H +D,\\
f(x_{n+1})\,&\in&\, f(x_{n+2}) + \gamma d(x_{n+1}, x_{n+2}) H +D,\\
\cdots\,&\  &\, \cdots \  \  \  \  \   \  \  \   \cdots  \  \  \  \   \   \  \   \   \cdots\\
f(x_{n+k-1})\,&\in&\, f(x_{n+k}) + \gamma d(x_{n+k-1}, x_{n+k}) H +D.
\end{eqnarray*}
Thus, there exist $h_{n+1},\, h_{n+2},\, \cdots,\, h_{n+k}\in H$ such that
\begin{eqnarray*}
f(x_n)\,&\in&\, f(x_{n+1}) + \gamma d(x_n, x_{n+1}) h_{n+1} +D,\\
f(x_{n+1})\,&\in&\, f(x_{n+2}) + \gamma d(x_{n+1}, x_{n+2}) h_{n+2} +D,\\
\cdots\,&\  &\, \cdots \  \  \  \  \   \  \  \   \cdots  \  \  \  \   \   \  \   \   \cdots\\
f(x_{n+k-1})\,&\in&\, f(x_{n+k}) + \gamma d(x_{n+k-1}, x_{n+k}) h_{n+k} +D.
\end{eqnarray*}
By adding the two sides of the above $k$ belonging relations, we have
\begin{eqnarray*}
f(x_n)\,&\in&\, f(x_{n+k}) + \gamma \sum_{i=n}^{n+k-1} d(x_i, x_{i+1}) h_{i+1} + D\\
&=&\, f(x_{n+k}) + \gamma \sum_{i=n+1}^{n+k} d(x_{i-1}, x_i) h_i +D.
\end{eqnarray*}
From this,
$$f(x_{n+k})\,\in\, f(x_n) - \gamma \sum_{i=n+1}^{n+k} d(x_{i-1}, x_i)h_i -D.\eqno{(27)}$$
Since $f$ is submonotone, we have $f(u) \leq_D  f(x_{n+k})$.
Combining this with (27), we have
$$f(u) \,\in\,  f(x_{n+k})-D \,\subset\, f(x_n) -\gamma \sum_{i=n+1}^{n+k} d(x_{i-1}, x_i) h_i-D. \eqno{(28)}$$
Remarking that  $H$ is $\sigma$-convex, we conclude that there exists $h_n^{\prime}\in H$ such that
$$\sum_{i=n+1}^{n+k} d(x_{i-1}, x_i) h_i\,\rightarrow\, \left(\sum_{i=n+1}^{\infty} d(x_{i-1}, x_i)\right) h_n^{\prime}\ \ \  (k\rightarrow\infty). \eqno{(29)}$$
Since $D$ is closed, by (28) and (29) we have
$$f(u)\,\in\, f(x_n) -\gamma \left(\sum_{i=n+1}^{\infty} d(x_{i-1}, x_i)\right) h_n^{\prime} -D$$
and so
$$f(x_n)\,\in\, f(u) + \gamma \left(\sum_{i=n+1}^{\infty} d(x_{i-1}, x_i)\right) h_n^{\prime} +D.\eqno{(30)}$$
On the other hand,
\begin{eqnarray*}
d(x_n, u)\, & =&\, \lim\limits_{k\rightarrow\infty} d(x_n, x_{n+k})\\
&\leq&\, \lim\limits_{k\rightarrow \infty}(d(x_n, x_{n+1}) + d(x_{n+1}, x_{n+2}) +
\cdots + d(x_{n+k-1}, x_{n+k}))\\ &=&\, \sum\limits_{i=n+1}^{\infty} d(x_{i-1}, x_i). \hspace{7.6cm}(31)
\end{eqnarray*}
By (30) and (31) we have
\begin{eqnarray*}
f(x_n) \,&\in&\, f(u)  + \gamma \left(\sum_{i=n+1}^{\infty}d(x_{i-1}, x_i)\right) h_n^{\prime} +D\\
&\subset&\, f(u) + \gamma d(x_n, u) h_n^{\prime}+D\\
&\subset&\, f(u)  +\gamma d(x_n, u) H +D.
\end{eqnarray*}
Thus, $u\in S(x_n)$ and assumption (C) is satisfied. Now, applying Theorem 2.1, there exists $\hat{x}\in X$ such that $\hat{x}\in S(x_0)$ and $S(\hat{x})  = \{\hat{x}\}$.
From this, we can easily show that $\hat{x}$ satisfies (a), (b) and (c).\hfill\framebox[2mm]{}\\

As is well-known,  for locally convex spaces, there are various notions of completeness. The weakest one seems to be local completeness (see [35, 36, 45]). A locally convex space  $Y$ is locally complete iff it is  $l^1$-complete, i.e., for each bounded sequence  $(b_n)\subset Y$ and each  $(\lambda_n) \subset l^1$, the series  $\sum_{n=1}^{\infty} \lambda_n b_n$  converges in $Y$.  Thus, if $Y$ is a locally complete locally convex space and $H$ is a locally closed, bounded convex set (or, $H$ is a cs-closed, bounded convex set), then we can show that $H$ is a $\sigma$-convex set. Concerning local completeness and local closedness, please refer to [35, Chapter 5] and [36, 37, 45]. Concerning cs-completeness and cs-closedness, please refer to [49].

As we have seen, the assumption in [40, Theorem 6.8] that there exists $\xi\in  D^+\cap H^{+s}$  such that $(f(X)-f(x_0))\cap (-(\cup_{\lambda >0}\lambda H +D))$ is $\xi$-lower bounded
has been replaced by here one that there exists $\epsilon >0$  such that
$(f(x_0) -\epsilon H-D) \cap  f(X)\,=\,\emptyset.$ We shall see that the latter is strictly weaker than the former from the following Proposition 5.2 and Example 5.3.  Hence, Theorem 5.1 improves [40, Theorem 6.8], and
also improves [7, Theorem 4.1] and [47, Theorem 6.2].\\

{PROPOSITION 5.2.}   \  {\sl Assume that there exists $\xi\in D^+\cap H^{+s}$ such that
$(F(X)-f(x_0))\cap(-(\cup_{\lambda >0}\lambda H +D))$ is $\xi$-lower bounded. Then, there exists $\epsilon >0$ such that $(f(x_0) -\epsilon H -D) \cap f(X) \,=\,\emptyset$.}\\

{\sl Proof.} \  If not, for every $n\in {\bf N}$,
$$(f(x_0)-n H -D) \cap f(X)\,\not=\,\emptyset.$$
From this,
$$(f(X)-f(x_0))\cap (-nH -D)\,\not=\,\emptyset,\ \ \forall n.$$
For each $n$,  there exists $z_n \,\in\, f(X)- f(x_0)$  such that
$$z_n\,\in\, -nH -D.\eqno{(32)}$$
Clearly,
$$z_n\,\in\,(f(X)-f(x_0)) \cap (-(\cup_{\lambda >0}\lambda H +D)).\eqno{(33)}$$
Since $\xi\in D^+\cap H^{+s}$, we have $\xi(d)\geq 0,\ \forall d\in D$ and $\alpha:=\inf\{\xi(h):\, h\in H\} > 0$.
Combining this with (32), we have
$$\xi(z_n)\,\leq\,-n\alpha,\ \ \forall n.$$
This with (33) contradicts the assumption that $(f(X)- f(x_0)) \cap (-(\cup_{\lambda >0} \lambda H+D))$ is $\xi$-lower bounded.\hfill\framebox[2mm]{}\\

The following example shows that there is  such a vector-valued map $f:\, X\rightarrow Y$ and $x_0\in X$ such that there exists $\epsilon >0$ such that $(f(X)-f(x_0))\cap (-\epsilon H-D)\,=\,\emptyset$, but for every  $\xi\in D^+\cap H^{+s}$,
$(f(X)- f(x_0))\cap (-(\cup_{\lambda >0}\lambda H +D))$ is not $\xi$-lower bounded.\\

{Example  5.3.}\  Let $X$ be ${\bf R}$ with the usual metric, i.e., $d(x, x^{\prime}) = |x-x^{\prime}|, \ x, x^{\prime}\in {\bf R}$, let $Y$ be ${\bf R}^2$ with the usual topology and with the partial order generated by the closed convex pointed cone $D \,=\,\{(y_1, y_2) \in {\bf R}^2:\, y_1\geq 0,\, y_2\geq 0\}$, and let $H\subset D\backslash\{0\}$ be a singleton $H=\{k_0\}$, where $k_0 = (1,1)\in D\subset {\bf R}^2$.
For any $\xi\in D^+\cap H^{+s} = D^+\cap \{k_0\}^{+s}$, there exists a unique $(\alpha, \beta) \in {\bf R}^2$ such that
$$\xi(y)\,=\,\alpha y_1 +\beta y_2,\ \  \forall y = (y_1, y_2)\in Y= {\bf R}^2.$$ From $\xi\in D^+\cap \{k_0\}^{+s}$, we conclude that $\alpha \geq 0,\, \beta \geq 0$ and at least one of $\alpha$ and $\beta$ is strictly greater than $0$, i.e., $\alpha >0,\, \beta \geq 0$ or $\alpha \geq 0,\, \beta >0$.
Let $f:\, X={\bf R}\,\rightarrow\, Y= {\bf R}^2$ be defined as follows:
 $$
f(x) = \left\{
\begin{array}{cc}
(-x, -1),\ \ \ &  \ \  {\rm if}  \ \ \  x>0;\\
(0, 0), \ \ \ & \ \  {\rm if}  \ \ \   x=0;\\
(-1, x), \ \ \ & \ \ {\rm if} \ \ \  x<0.
\end{array}
\right.
$$
Put $x_0:=0$ and $\epsilon =2$.
Then
\begin{eqnarray*}
f(X)-f(x_0)\,&=&\, \{(-x, -1): x>0\}\cup\{(-1, x): x<0\} \cup\{(0,0)\}\\
&=&\,\{(x, -1): x<0\}\cup\{(-1, x): x<0\} \cup\{(0,0)\}.
\end{eqnarray*}
Also,
\begin{eqnarray*}
-\epsilon H -D\,&=&\,-2k_0 -D\\
&=&\, -2 (1,1)  + \{(y_1, y_2) \in {\bf R}^2:\, y_1\leq 0,\, y_2\leq 0\}\\
&=&\,\{(y_1-2,\, y_2-2):\, y_1\leq 0, y_2\leq 0\}\\
&=&\,\{(y_1, y_2):\, y_1\leq -2, \, y_2\leq -2\}.
\end{eqnarray*}
Obviously,
$$(f(X) - f(x_0))\,\cap\,(-\epsilon H -D)\,=\,\emptyset.$$
On the other hand,
\begin{eqnarray*}
&\  &\, (f(X) - f(x_0))\cap (-(\cup_{\lambda >0}\lambda H +D))\\
&=&\, (f(X)- f(x_0))\cap (-(\cup_{\lambda >0} \lambda k_0 +D))\\
&=&\, \left(\{(x, -1): x<0\} \cup\{(-1, x): x<0\} \cup\{(0,0)\}\right)\\
&\  &\, \cap\, \{(y_1, y_2)\in {\bf R}^2:\, y_1 < 0, y_2 <0\}\\
&=&\,\{(x,-1):\, x<0\} \cup \{(-1, x):\, x<0\}.
\end{eqnarray*}
For any $\xi\in D^+\cap \{k_0\}^{+s}$, there exists a unique $(\alpha, \beta) \in {\bf R}^2$ such that
$\xi(y) = \alpha y_1 +\beta y_2,\ \forall y=(y_1, y_2)\in Y = {\bf R}^2$, where $\alpha \geq 0,\ \beta \geq 0$ and at least
one of $\alpha,\ \beta$ is strictly greater than $0$.
Thus,
\begin{eqnarray*}
&\ &\, \xi\circ \left((f(X)-f(x_0))\cap (-(\cup_{\lambda >0}\lambda H +D))\right)\\
&=&\, \{\alpha x-\beta:\, x<0\} \cup \{-\alpha + \beta x:\, x<0\},
\end{eqnarray*}
where $\alpha \geq 0$ and $\beta \geq 0$.
If $\alpha >0$, then $\{\alpha x-\beta:\, x<0\}$ is not lower bounded.
If $\beta >0$,  then $\{-\alpha +\beta x:\, x <0\}$ is not lower bounded.
Hence, for any $\xi\in D^+\cap H^{+s}$,
$$(f(X)-f(x_0)) \cap (-(\cup_{\lambda >0}\lambda H +D))$$
is not $\xi$-lower bounded.\\

\newpage
\noindent{\bf References} \vskip 10pt
\begin{description}
\small
\item{[1]}  M. Adan, V. Novo, Proper efficiency in vector optimization
on real linear spaces, J. Optim. Theory Appl., 121 (2004), 515-540.

\item{[2]}  Y. Araya, Ekeland's variational principle and its
equivalent theorems in vector optimization, J. Math. Anal. Appl.,
 346 (2008), 9-16.

\item{[3]} J. P. Aubin,  H. Frankowska, Set-Valued Analysis,
Birkh\"{a}user, Boston,  1990.

\item{[4]} T. Q. Bao, B. S. Mordukhovich, Variational principles for
set-valued mappings with applications to multiobjective
optimization, Control Cybern., 36 (2007), 531-562.

\item{[5]} T. Q. Bao, B. S. Mordukhovich, Relative Pareto minimizers
for multiobjective problems: existence and optimality conditions,
Math. Program, Ser.A, 122 (2010), 301-347.

\item{[6]}  E. M. Bednarczuk,  M. J. Przybyla, The vector-valued
variational principle in Banach spaces ordered by cones with
nonempty interiors, SIAM J. Optim. 18 (2007), 907-913.

\item{[7]}  E. M. Bednarczk, D. Zagrodny, Vector variational
principle, Arch. Math. (Basel), 93 (2009), 577-586.

\item{[8]}  Y. Chen,  Y. J. Cho,  L. Yang, Note on the results with lower semi-continuity, Bull. Korean Math. Soc.,  39 (2002),  535-541.

\item{[9]} G. Y. Chen, X. X. Huang, A unified approach to the existing three
types of variational principle for vector valued functions, Math.
Methods Oper. Res., 48 (1998), 349-357.

\item{[10]}  G. Y. Chen, X. X. Huang, X. G. Yang, Vector Optimization, Set-Valued
and Variational Analysis, Springer-Verlag, Berlin, 2005.

\item{[11]} D. Dentcheva, S. Helbig, On variational principles,
level sets, well-posedness, and $\epsilon$-solutions in vector
optimization, J. Optim. Theory Appl., 89 (1996), 325-349.

\item{[12]}  W. S. Du, On some nonlinear problems induced by an abstract maximal element principle, J. Math. Anal. Appl., 347 (2008), 391-399.

\item{[13]}  I. Ekeland, Sur les prob\`{e}mes variationnels, C. R. Acad. Sci.
Paris 275 (1972), 1057-1059.

\item{[14]}  I. Ekeland,  On the variational principle, J. Math.
Anal. Appl. 47 (1974), 324-353.

\item{[15]} I. Ekeland, Nonconvex minimization problems, Bull. Amer.
Math. Soc. (N.S.) 1 (1979), 443-474.

\item{[16]}  C. Finet, L. Quarta, C. Troestler,  Vector-valued
variational principles, Nonlinear Anal. 52 (2003), 197-218.

\item{[17]}  F. Flores-Baz\'{a}n, C. Guti\'{e}rrez,  V. Novo,  A
Br\'{e}zis-Browder principle on partially ordered spaces and related
ordering theorems, J. Math. Anal. Appl., 375 (2011), 245-260.

\item{[18]} Chr. Gerstewitz (Tammer), Nichtkonvexe Dualit\"{a}t in
der Vektoroptimierung, Wiss. Z. TH Leuna-Merseburg 25 (1983),
357-364.

\item {[19]}  A. G\"{o}pfert, H. Riahi,  Chr. Tammer,  C.
Z$\breve{a}$linescu,  Variational Methods in Partially Ordered
Spaces, Springer-Verlag, New York, 2003.

 \item{[20]} A. G\"{o}pfert,  C. Tammer and C. Z$\breve{a}$linescu,  On the
vectorial Ekeland's variational principle and minimal point theorems
in product spaces, Nonlinear Anal. 39 (2000),  909-922.

\item{[21]} C. Guti\'{e}rrez, B. Jim\'{e}nez, V. Novo, On
approximate efficiency in multiobjective programming, Math. Methods
Oper. Res., 64(2006),  165-185.

\item{[22]} C. Guti\'{e}rrez, B. Jim\'{e}nez, V. Novo,  A unified
approach and optimality conditions for approximate solutions of
vector optimization problems, SIAM J. Optim., 17(2006),  688-710.

\item{[23]} C. Guti\'{e}rrez, B. Jim\'{e}nez, V. Novo, A set-valued
Ekeland's variational principle in vector  optimization, SIAM J.
Control. Optim., 47 (2008), 883-903.

\item{[24]} T. X. D. Ha, Some variants of the Ekeland  variational
principle for a set-valued map, J. Optim. Theory Appl., 124 (2005),
187-206.

\item{[25]}  A. H. Hamel, Equivalents to Ekeland's variational
principle in uniform spaces, Nonlinear Anal. 62 (2005), 913-924.

\item{[26]}  F. He, J. H. Qiu, Sequentially lower complete spaces
and Ekeland's variational principle,  Acta Math. Sinica. English
Series, 31 (2015), 1289-1302.

\item{[27]}  J. Horv\'{a}th, Topological Vector Spaces and Distributions, vol. 1,
Addison-Wesley, Reading, MA, 1966.

\item{[28]} G. Isac, The Ekeland's principle and the
Pareto $\epsilon$-efficiency, In: M. Tamiz (ed.) Multi-Objective
Programming and Goal Programming: Theories and Applications, Lecture
Notes in Econom. and Math. Systems, vol. 432, Springer-Verlag,
Berlin, 1996, 148-163.

\item{[29]} P. Q. Khanh,  D. N. Quy, Versions of Ekeland's variational principle involving set perturbations, J. Glob. Optim.,  57 (2013), 951-968.

\item{[30]}  J. L. Kelley, I. Namioka, et al., Linear Topological
Spaces, Van\ Nostrand, Princeton, 1963.

\item{[31]}  G.  K\"{o}the, , Topological Vector Spaces I, Springer-Verlag,
Berlin, 1969.

\item{[32]}  C. G. Liu, K. F. Ng, Ekeland's variational principle for
set-valued functions, SIAM J. Optim., 21 (2011), 41-56.

\item{[33]}  A. B. N\'{e}meth, A nonconvex vector minimization problem,
Nonlinear Anal. 10 (1986),  669-678.

\item{[34]}  D. Pallaschke, S. Rolewicz, Foundations of Mathematical
Optimization, Math. Appl. 388, Kluwer, Dordrecht, 1997.

\item{[35]}  P. P\'{e}rez Carreras, J. Bonet, Barrelled Locally Convex Spaces,
North-Holland, Amsterdam, 1987.

\item{[36]} J. H. Qiu, Local completeness and dual local quasi-completeness,
Proc. Amer. Math. Soc. 129 (2001) 1419-1425.

\item{[37]} J. H. Qiu, Local completeness and drop theorem, J. Math.
Anal. Appl. 266 (2002), 288-297.

\item{[38]}  J. H. Qiu, A generalized Ekeland vector variational
principle and its applications in optimization, Nonlinear Anal., 71
(2009), 4705-4717.

\item{[39]} J. H. Qiu, On Ha's version of set-valued Ekeland's
variational principle, Acta Math. Sinica, English Series, 28 (2012),
717-726.

\item{[40]} J. H. Qiu, Set-valued quasi-metrics and a general
Ekeland's variational principle in vector optimization, SIAM J.
Control  Optim., 51 (2013), 1350-1371.

\item{[41]} J. H. Qiu, The domination property for efficiency and
Bishop-Phelps theorem in locally convex spaces, J. Math. Anal.
Appl., 402 (2013), 133-146.

\item{[42]} J. H. Qiu, A pre-order principle and set-valued  Ekeland variational principle, J. Math. Anal. Appl., 419 (2014), 904-937.

\item{[43]} J. H. Qiu, F. He, A  general vectorial Ekeland's
variational principle with a p-distance, Acta Math. Sinica, 29
(2013), 1655-1678.

\item{[44]} J. H. Qiu, B. Li, F. He, Vectorial Ekeland's variational
principle with a w-distance and its equivalent theorems, Acta Math.
Sci., 32B (2012), 2221-2236.

\item{[45]} S. A. Saxon, L. M. S\'{a}nchez Ruiz, Dual local
completeness, Proc. Amer. Math. Soc., 125 (1997), 1063-1070.

\item{[46]} C. Tammer, A generalization of Ekeland's variational principle,
Optimization,  25 (1992)  129-141.

\item{[47]} C. Tammer, C. Z$\breve{a}$linescu, Vector variational
principle for set-valued functions, Optimization, 60 (2011),
839-857.

\item{[48]}  A. Wilansky, Modern Methods in Topological Vector Spaces,
McGraw-Hill, New York, 1978.

\item{[49]} C. Z$\breve{a}$linescu, Convex Analysis  in   General
Vector Spaces, World Sci., Singapore, 2002.

\item{[50]} J. Zhu, C. K. Zhong, Y. J. Cho, A generalized variational principle
and vector optimization, J. Optim. Theory Appl., 106 (2000) 201-218.

\end{description}
\end{document}